\documentclass[titlepage,11pt]{article}
\usepackage{latexsym}
\usepackage{amsfonts}
\usepackage{amsmath}
\newcommand\blackslug{\hbox{\hskip 1pt \vrule width 4pt height 8pt depth 1.5pt
        \hskip 1pt}}
\newcommand\bbox{\hfill \quad \blackslug \bigbreak}

\def\ll{,\ldots,}

\title{Small families under subdivision}
\author{Maria Chudnovsky\thanks{This material is based upon work supported in part by the U. S. Army
Research Office under grant   number W911NF-16-1-0404, and supported by  NSF grant DMS 1763817.}\\
Princeton University, Princeton, NJ 08544, USA
\\
\\
Martin Loebl\thanks{Partially supported by the H2020-MSCA-RISE project CoSP- GA No. 823748.}\\
Charles University, Prague, Czech Republic 
\\
\\
Paul Seymour\thanks{Partially supported by NSF grant  DMS-1800053 and AFOSR grant A9550-19-1-0187.}\\
Princeton University, Princeton, NJ 08544, USA}

\date{March 23, 2019; revised \today}

\newtheorem{thm}{}[section]

\newcommand{\Proof}{\noindent{\bf Proof.}\ \ }

\begin{document}
\maketitle
\begin{abstract}
Let $H$ be a graph with maximum degree
$d$, and let $d'\ge 0$. We show that for some $c>0$ depending on $H,d'$, and all integers $n\ge 0$, there are at most $c^n$ unlabelled simple $d$-connected 
$n$-vertex graphs with maximum degree at most $d'$ that do not contain
$H$ as a subdivision. On the other hand, the number of unlabelled simple $(d-1)$-connected $n$-vertex graphs with minimum degree $d$ and
maximum degree at most
$d+1$ that do not contain
$K_{d+1}$ as a subdivision is superexponential in $n$.
\end{abstract}

\section{Introduction}
All graphs in this paper are finite and have no loops or parallel edges.
 We say $H$ is a {\em minor} of $G$ if $H$ can be obtained from a  subgraph of $G$ by contracting edges.
A {\em subdivision} of a graph $H$ is a graph $J$ such that $H$ can be obtained from $J$ by repeatedly contracting an edge
incident with a vertex of degree two; and we say $G$ {\em contains a subdivision of $H$} if some subgraph of $G$ is a subdivision of $H$.

It is known~\cite{amini,norin} that for every graph $H$, the number of unlabelled $n$-vertex graphs
not containing $H$ as a minor is at most $c^n$, for some constant $c$ depending on $H$, as we discuss in the next section.
Here we investigate the number of 
$n$-vertex graphs that do not contain a subdivision of $H$.

When is it true that the number of (unlabelled)
$n$-vertex graphs that do not contain $H$ as a subdivision is only exponential in $n$? When $H$ has maximum degree at most three, 
this is true, since in that case containing a subdivision of $H$ is equivalent to containing $H$ as a minor. But when $H$ has maximum
degree at least four, it is false: in fact there are superexponentially many graphs of maximum degree three that do not contain $H$.
So, this question is not very interesting. Let us ask instead, if $H$ has maximum degree $d$, when is it true that the number
of $n$-vertex graphs with minimum degree at least $d$ not containing a subdivision of $H$ is at most exponential?

This is still not true in general. For instance, let $G$ be four-connected, with a vertex $v$ 
such that $G\setminus v$ is a cubic graph.
Then $G$ contains no subdivision of the octahedron graph (the 4-regular graph on six vertices); 
and yet there are superexponentially many such graphs $G$, since there
are superexponentially many cubic graphs that become 4-connected if we 
add one extra vertex adjacent to all other vertices.  For this reason we will impose an upper bound on degree.

This still is not enough: we shall show in section~\ref{sec:gadget} that 
\begin{thm}\label{firstgadget}
For all $d\ge 5$, and infinitely many $n>0$,
there are superexponentially many $d$-regular $(d-1)$-connected $n$-vertex graphs
that do not contain a subdivision of $K_{d+1}$. 
\end{thm}
(This is not true for $d=4$.)
So let us assume $d$-connectivity as well. Then we have a positive result:
\begin{thm}\label{mainthm}
Let $d,d'\ge 0$ and let $H$ have maximum degree at most $d$. Then there exists $c>0$ such that there are at most $c^n$ unlabelled 
$d$-connected 
$n$-vertex graphs with maximum degree at most $d'$ that do not contain a subdivision of $H$.
\end{thm}

This is an immediate consequence of \ref{smallfam2} below, and the following:

\begin{thm}\label{subtominor}
Let $d,d'\ge 0$, and let $H$ be a graph with maximum degree at most $d$. 
Then there exists $t\ge 0$ such that for every graph $G$, if $G$ is $d$-connected 
and has maximum degree at most $d'$, and $G$ contains $K_t$ as a minor, then $G$ contains a subdivision of $H$.
\end{thm}

Again, this becomes false if we omit the upper bound $d'$ on degree; the same example (a cubic graph with an extra vertex
adjacent to all others) can have arbitrarily large complete graphs as minors, and this generalizes to other values of $d$.

As far as we know, the result \ref{subtominor} is new, despite the extensive research that has already been carried out in this area. 
Chun-Hung Liu (private communication) tells us that it can be deduced from theorem 6.8
in~\cite{liu}, but we will give a proof not using such heavy machinery. A related result was proved by Irene Muzi in her 
thesis~\cite{muzi}:
\begin{thm}\label{muzi}
 Let
$H_1$ be a graph of maximum degree at most four, and let $H_2$ be a graph such that for some vertex $v$, 
$H_2\setminus v$ has maximum degree at most three.  Then there exists $t$ such that every
4-connected graph that has a $K_t$ minor contains a subdivision of one of $H_1,H_2$.
\end{thm}

As we said, the statement of \ref{firstgadget} is false when $d=4$, and is worth a closer look. For infinitely many $n$ there are 
superexponentially many $3$-connected $n$-vertex graphs
that do not contain a subdivision of $K_5$, in which every vertex has degree $4$ or $5$; but 
the number of such graphs that are $4$-regular graphs is only exponential. That follows from:
\begin{thm}\label{4subtominor}
Let $d\in \{4,5\}$ and $d'\ge 0$, and let $H$ be a graph with maximum degree at most $d$. 
Then there exists $t\ge 0$ such that for every graph $G$, if $G$ is $(d-1)$-connected and $d$-edge-connected
and has maximum degree at most $d'$, and $G$ contains $K_t$ as a minor, then $G$ contains a subdivision of $H$.
\end{thm}
This provides a substitute for \ref{subtominor} for $4$-regular $3$-edge-connected graphs, since for $d=4$, $d$-regular $(d-1)$-connected graphs are
$d$-edge-connected.

\section{Minor containment}

It was proved in \cite{smallfam} that
\begin{thm}\label{smallfam} 
For every graph $H$, there exists $c$ such that for all $n$, the number of labelled $n$-vertex graphs
not containing $H$ as a minor is at most $n!c^n$.
\end{thm}

A strengthening was found by Amini, Fomin and Saurabh~\cite{amini}, and the same result is implicit in the proof of
theorem 4 of~\cite{norin}:
\begin{thm}\label{smallfam2}
For every graph $H$, there exists $c$ such that for all $n$, the number of unlabelled $n$-vertex graphs
not containing $H$ as a minor is at most $c^n$.
\end{thm}

This paper is concerned with subdivision containment rather than minor containment, but let us take the opportunity
to give another proof of \ref{smallfam2}, simpler than those in print. Two (distinct) vertices $u,v$ of a graph $G$
are {\em twins} if they have the same neighbours in $V(G)\setminus \{u,v\}$. We use the following easy lemma of~\cite{smallfam}:
\begin{thm}\label{lowdegree}
For every graph $H$ there exist $c,d$ with $0<c\le 1$ and $d\ge 1$ such that for all $n\ge 2$ and for every $n$-vertex graph $G$ not 
containing $H$ as a minor, 
there is a subset $S\subseteq V(G)$ with $|S|\ge cn$, with the following properties. Every vertex in $S$ has degree at most $d$;
and for every
$v\in S$, there exists $u\in S$ with $u\ne v$, such that either 
$u,v$ are nonadjacent twins, or
$u,v$ are adjacent.
\end{thm}

\noindent{\bf Proof of \ref{smallfam2}.\ \ }
If $J$ is a graph, we say that a graph $G$ is a {\em $k$-growth} of $J$ if there are $2k$ distinct vertices $u_1,v_1\ll u_k, v_k$ of $G$, such that either
\begin{itemize}
\item for $1\le i\le k$, $u_i,v_i$
are nonadjacent and have the same neighbours in $G$, and $J$ is obtained from $G$ by deleting $v_1\ll v_k$; or
\item for $1\le i\le k$, $u_i,v_i$
are adjacent and both have degree at most $d$, and $J$ is obtained from $G$ by contracting the edges $u_iv_i$ for $1\le i\le k$
(and deleting any resultant parallel edges).
\end{itemize}
Let us call these $k$-growths of the {\em first} and {\em second kind} respectively. 
\\
\\
(1) {\em If $J$ is an $m$-vertex graph and $k\ge 0$, the number of distinct unlabelled $k$-growths of $J$ is at most $3^{2dk}2^m$.}
\\
\\
Certainly $J$
has at most $2^m$ $k$-growths of the first kind, since every such $k$-growth is determined by a knowledge of the set
$\{u_1\ll u_k\}\subseteq V(J)$, and there are at most $2^m$ such sets.
We claim that $J$ has at most $2^m3^{(2d-2)k}$ $k$-growths of the second kind. To see this, let $G$ be a $k$-growth of the second kind,
let $u_i,v_i\in V(G)\;(1\le i\le k)$ be as in the definition, and let $w_1\ll w_k$
be the vertices of $J$ formed by contracting the edges $u_1v_1\ll u_kv_k$. Thus each $w_i$ has degree at most $2d-2$. There are
at most $2^m$ choices for the set $\{w_1\ll w_k\}$. For each such choice $\{w_1\ll w_k\}$, to construct $G$ 
we replace each $w_i$ by the edge $u_iv_i$, and for each edge $w_ix$ we replace it by either $u_ix$ or $v_ix$, or both.
Let us follow this process one edge at a time. We have three choices for each edge $w_ix$, and
there are only at most $(2d-2)k$ such edges since $w_1\ll w_k$ have degree at most $2d-2$ (and we may assume this property is preserved
as we do the uncontractions one by one). Thus, each choice of $w_1\ll w_k$ gives rise to at most $3^{(2d-2)k}$
$k$-growths of the second kind; and hence in total there are at most $3^{(2d-2)k}2^m$ $k$-growths of the second kind.
Since there are only $2^m$ $k$-growths of the first kind, we deduce that altogther there are at most 
$$3^{(2d-2)k}2^m+2^m\le 3^{2dk}2^m$$
$k$-growths of $J$. This proves (1).

\bigskip

Given $H$, let $c,d$ be as in \ref{lowdegree}. 
\\
\\
(2) {\em If $G$ is an $n$-vertex graph with no $H$ minor, then $G$ is a $k$-growth of some graph
where $k=\lceil cn/(2d+3)\rceil$.}
\\
\\
Let $G$ be a $n$-vertex graph with no $H$ minor, and let $S$ be as 
in \ref{lowdegree}. Let $S_2$ be the set of vertices in $S$ with a neighbour in $S$, and let $S_1=S\setminus S_2$.
For each $v\in S_1$ there exists $u\in S$ such that $u,v$ are nonadjacent twins; and consequently $u\in S_1$. Let us partition $S_1$
into pairs of twins, as far as possible; then we obtain at least $|S_1|/3$ pairs. Consequently $G$ is a $k_1$-growth
of the first kind, where $k_1\ge |S_1|/3$. The subgraph induced on $S_2$ has minimum degree at least one, and maximum
degree at most $d$, and so has a matching of cardinality at least $|S_2|/(2d)$, as is easily seen.
Hence $G$ is a $k_2$-growth of the second kind, where $k_2\ge |S_2|/(2d)$. Now 
$$cn\le |S_1|+|S_2|\le 3k_1+2dk_2\le (2d+3)\max(k_1,k_2)$$
and so $\max(k_1,k_2)\ge k$. This proves (2).

\bigskip

For $n\ge 0$, let there be $f(n)$ unlabelled $n$-vertex graphs not containing $H$ as a minor.
Define $b=3^{2d}2^{(2d+3)/c-1}$. Thus $b\ge 1$.
We prove by induction on $n$ that $f(n)\le b^n$.  This is true if $n\le 1$, so we assume that $n\ge 2$. Let $k$ be as in (2).
By (1) and (2), $f(n)\le f(n-k)3^{2dk}2^{n-k}$. From the inductive hypothesis, $f(n-k)\le b^{n-k}$; and so
$$f(n)b^{-n}\le b^{-k}3^{2dk}2^{n-k}.$$
It remains to show that $b^{-k}3^{2dk}2^{n-k}\le 1$, that is, that $3^{2d}2^{n/k-1}\le b$. But $n/k\le (2d+3)/c$ from the definition
of $k$, and so $3^{2d}2^{n/k-1}\le 3^{2d}2^{(2d+3)/c-1}=b$. This proves \ref{smallfam2}.~\bbox

\section{A construction}\label{sec:gadget}

In this section we prove \ref{firstgadget}.
We need 
\begin{thm}\label{manygraphs}
For all $b\ge 0$, and for infinitely many $m>0$, there are more than $b^m$ different unlabelled
graphs with $m$ vertices that are $d$-regular and
$d$-connected.
\end{thm}

In fact Bender and Canfield~\cite{bender}, and independently Wormald~\cite{wormald}, proved the following, which for $d\ge 3$ is much stronger than \ref{manygraphs}:
\begin{thm}\label{bender}
For all integers $d\ge 0$ and $m\ge 0$ with $dm$ even, the number of labelled $d$-regular graphs with $m$ vertices is approximately
$$\sqrt{2}e^{1-d^2/4}\left(\frac{d^dm^d}{e^d(d!)^2}\right)^{m/2},$$
where ``$f(m)$ is approximately $g(m)$'' means $f(m)/g(m)\rightarrow 1$ as $m\rightarrow\infty$.
\end{thm}
(In fact these papers estimated the number of labelled $d$-regular graphs; to count just the $d$-connected ones, a result of  
\L uczak~\cite{luczak} shows that for fixed $d\ge 3$, as $m\rightarrow\infty$ almost all labelled $d$-regular graphs 
with $m$ vertices are $d$-connected.) To see that \ref{bender} implies \ref{manygraphs}, observe that for any constant $b$,
 the expression in \ref{bender}
is more than $m!b^m$ for $m$ sufficiently large.

The next result is a lemma that will help to show that the graphs we construct do not contain a subdivision of $K_{d+1}$.
\begin{thm}\label{littlebag}
Let $G$ be a $d$-regular graph, and let $Z\subseteq V(G)$ with $|Z|=2d-2$. Let $X\subseteq Z$ with $|X|\le d-1$,
such that at least three vertices in $X$ have exactly one neighbour in 
$V(G)\setminus Z$, and each vertex in $Z\setminus X$ has no neighbours in $V(G)\setminus Z$. Let $H$ be a subgraph of $G$
that is a subdivision of $K_{d+1}$, and let $U$ be the set of vertices that have degree $d$ in $H$.
Then every vertex in $Z\cap U$ belongs to $X$
and has at least two neighbours in $V(G)\setminus Z$.
\end{thm}
\Proof
Let $U=\{u_1\ll u_{d+1}\}$. Then evidently:
\\
\\
(1) {\em For $1\le i<j\le d+1$ there are $d$ paths of $H$ between $u_i$ and $u_j$, pairwise internally disjoint.}
\\
\\
(2) {\em Every vertex of $G$ adjacent in $G$ to at least three members of $U$ belongs to $U$.}
\\
\\
For $1\le i \le d+1$, every edge of $G$ incident with $u_i$ is an edge of $H$, since $G$ is $d$-regular and 
$u_i$ has degree $d$ in $H$. So every vertex of $G$ with a neighbour in $U$ belongs to $V(H)$.
Suppose that $v\in V(G)$
is adjacent in $G$ to at least three members of $U$. Then $v\in V(H)$, and has degree at least three in $H$; 
and since $H$ is a subdivision
of $K_{d+1}$, it follows that $v\in U$. This proves (2).

\bigskip
Suppose first that $U\cap (Z\setminus X)\ne \emptyset$, and let $u_1\in (Z\cap U)\setminus X$. 
Since $X$ is a cutset of cardinality at most $d-1$, (1) implies that $U\subseteq Z$. Let $X'$ be the set of all vertices in $X$
with at least two neighbours in $V(G)\setminus Z$.
For each $v\in Z\setminus X'$, since $v$ has at most one neighbour in $V(G)\setminus Z$, it has at least $d-1$ neighbours in $Z$, 
and since $|Z\setminus U|=d-3$, it
follows that either $v\in U$ or $v$ has at least three neighbours in $U$; and from (2), it follows that $v\in U$. But then
$U$ includes $Z\setminus X$ and all the (at least three) vertices of $X\setminus X'$, which is impossible since $|Z\setminus X|\ge d-1$
and $|U|=d+1$. This proves that $U\cap Z\subseteq X$, and we claim that $U\cap Z\subseteq X'$. Suppose not, and
let $u\in (U\cap Z)\setminus X'$, and let $y$ be the unique neighbour of $u\in V(G)\setminus Z$. Since $|U|=d+1$ and $|X|\le d-1$,
there are at least two vertices of $U$ not in $Z$, and so there is one of them, say $u'$, that is different from $y$.
Then $(X\setminus \{u\})\cup \{y\}$ is a cutset of cardinality
at most $d-1$ separating $u$ and $u'$, contrary to (1). This proves \ref{littlebag}.~\bbox

Now we prove \ref{firstgadget}, which we restate:
\begin{thm}\label{gadget}
Let $d\ge 5$ be an integer; then for all $c>0$, and infinitely many values of $n>0$, there are more than $c^n$
$d$-regular $(d-1)$-connected graphs $G$ with $n$ vertices that do not contain a subdivision of $K_{d+1}$.
\end{thm}
\Proof
The proof breaks into two cases, depending whether $d$ is odd or even. The odd case is easier, so we begin with that.
Take a complete bipartite graph with bipartition $(X,Y)$, where $|X|=|Y|=d-1$. Partition $Y$ into $(d-1)/2$ pairs,
and add an edge joining each pair, forming a graph $R$ say. Thus every vertex in $Y$ has degree $d$, and every vertex in $X$
has degree $d-1$.

Let $n=(2d-2)m$, where $m\ge 2$ is an integer.
Let $D$ be a $(d-1)$-regular, $(d-1)$-connected graph with $m$ vertices $v_1\ll v_m$. Take the disjoint union of 
$m$ copies of $R$, say $R_1\ll R_m$, and for $1\le i\le n$ let $X_i\subseteq V(R_i)$ be the set of vertices of $R_i$
with degree $d-1$ in $R_i$, and $Y_i=V(R_i)\setminus X_i$.
For each edge $v_iv_j$ of $D$, add an edge between $X_i$ and $X_j$, so that these new edges form a matching (this is possible,
since each $v_i$ has degree $d-1$ in $D$, and each set $X_i$ has cardinality $d-1$). Let the graph we form by this process be $G_D$.
It is easy to see that $G_D$ is $(d-1)$-connected and $d$-regular, and has $n$ vertices. We claim that $G_D$ contains 
no subdivision of $K_{d+1}$. Suppose that it does, and $H\subseteq G_D$ is a subdivision of $K_{d+1}$. Let $U$
be the vertices of $H$ with degree $d$. By \ref{littlebag} applied to $Z=V(R_i)$, it follows that $U\cap V(R_i)=\emptyset$,
for $1\le i\le m$, which is impossible.
This proves that $G_D$
contains no subdivision of $K_{d+1}$ (in the case when $d$ is odd).

We observe that if $D,D'$ are nonisomorphic graphs, both $(d-1)$-regular and $(d-1)$-connected, then $G_D$ and $G_{D'}$ are 
nonisomorphic; because
with notation as before, every subgraph of $G_D$ isomorphic to $R$ is one of the graphs $R_1\ll R_m$, and so $D$
can be obtained from $G_D$ by contracting all edges of every $R$-subgraph of $G_D$. 

Now let $c>0$, and let $b=c^{2d-2}$. By \ref{manygraphs}, for infinitely many $m>0$
there are more than $b^m$ graphs $D$ on $m$ vertices that are $(d-1)$-regular and $(d-1)$-connected; and 
all the corresponding graphs $G_D$ are distinct. Hence, since each such $G_D$ has $(2d-2)m=n$ vertices,
for infinitely many $n>0$ there are more than $b^m=c^n$ $d$-regular $(d-1)$-connected graphs with $n$ vertices,
not containing a subdivision of $K_{d+1}$. This completes the proof when $d$ is odd.

Now we turn to the case when $d$ is even. Define a graph $R$ as follows. Take six pairwise disjoint sets 
$A,B,C, C', B',A'$ of cardinalities $d/2,d-1,d/2-1,d/2-1,d-1,d/2$ respectively, and
choose $s,u\in B, t\in C, t'\in C'$ and $s',u'\in B'$. Make every vertex in $A\cup C$
adjacent to every vertex in $B$, except the pair $st$; and similarly make $A'\cup C'$ complete to $B'$ except for the pair $s't'$.
Add a perfect matching between $C,C'$, in which $t, t'$ are not adjacent (this is possible since $d/2-1\ge 2$);
and add one more edge $tt'$.
Pair up the vertices in $B\setminus \{u\}$, and join each pair with an edge; and add one more edge $su$. (Thus $s$ is incident with two
of these edges.)
Add edges within $B'$ similarly. This defines $R$. We see that
\begin{itemize}
\item every vertex in $A\cup A'$ has degree $d-1$, and all other vertices have degree $d$;
\item every edge of $R$ either has both ends in $A\cup B\cup C$, or both ends in $A'\cup B'\cup C'$,
or both  ends in $C\cup C'$.
\end{itemize}

Let $n=(4d-4)m$, where $m\ge 2$ is an integer.
Let $D$ be a $d$-regular, $d$-connected graph with $m$ vertices $v_1\ll v_m$. Take the disjoint union of
$m$ copies of $R$, say $R_1\ll R_m$, and for $1\le i\le m$
for $1\le i\le m$, let
$$A_i,B_i,C_i,C_i', B_i', A_i', s_i, t_i, u_i$$
correspond to
$$A,B,C,C', B', A', s, t, u$$
respectively.
For each edge $v_iv_j$ of $D$, add an edge between $X_i$ and $X_j$, so that these new edges form a matching (this is possible,
since each $v_i$ has degree $d$ in $D$, and each set $X_i$ has cardinality $d$.) Let the graph we form by this process be $G_D$.
It is easy to see that $G_D$ is $d$-regular, and has $n$ vertices. We claim that $G_D$ contains
no subdivision of $K_{d+1}$. Suppose that it does, and $H\subseteq G_D$ is a subdivision of $K_{d+1}$. Let $U$
be the vertices of $H$ with degree $d$. Choose some $i$ with $1\le i\le m$.
Since every vertex in $A_i$ has only one neighbour in $V(G)\setminus Z$
and $|A_i|=d/2\ge 3$, \ref{littlebag} applied to $Z=A_i\cup B_i\cup C_i$ tells us that $U\cap (A_i\cup B_i)=\emptyset$,
and similarly $U\cap (A_i'\cup B_i')=\emptyset$. Suppose that there exists $u\in U\cap C_i$. Since $A_i\cup C_i'$
is a cutset of $G$ of cardinality $d-1$, and every two vertices in $U$ are joined by $d$ internally disjoint paths, it follows that
$U\subseteq A_i\cup B_i\cup C_i\cup C_i'$, and hence $U\subseteq C_i\cup C_i'$. But this is impossible since $|C_i\cup C_i'|=d-2$.
Thus $U\cap C_i=\emptyset$, and similarly $U\cap C_i'=\emptyset$; and so $U\cap V(R_i)=\emptyset$. Since this holds for all $i$,
this is impossible. Consequently $G_D$ contains no subdivision of $K_{d+1}$.

For $1\le i\le m$, let $X_i=A_i\cup A_i'$, 
and $Y_i=V(R_i)\setminus X_i$.
Next we show that the graph $G_D$ is $(d-1)$-connected.
Suppose not; then there is a set $W\subseteq V(G_D)$
with $|W|\le d-2$ such that deleting $W$ from $G_D$ makes a graph with at least two components. Choose a partition $P,Q$
of $V(G_D)\setminus W$ with $P,Q$ nonempty, such that no vertex in $P$ has a neighbour in $Q$.
\\
\\
(1) {\em For $1\le i\le m$, not both $P,Q$ have nonempty intersection with $A_i\cup B_i\cup C_i$.}
\\
\\
Suppose they do, for $i = 1$ say. Since $|A_1\cup C_1|=|B_1|=d-1$ and $|W|\le d-2$, and $P\cup Q$ contains all vertices not in $W$,
it follows that $P\cup Q$ has nonempty intersection with 
both $A_1\cup C_1$ and $B_1$. If $P\cup Q$ contains a vertex of $A_1\cup C_1$ different from $s_1$, and contains a vertex of $B_1$
different from $t_1$, then they are adjacent, and since every vertex of $A_1\cup B_1\cup C_1$ is adjacent to one of these two
vertices, it follows that the subgraph induced on $(P\cup Q)\cap (A_1\cup B_1\cup C_1)$ is connected, a contradiction.
Thus either $W$ includes $A_1\cup C_1\setminus \{t_1\}$ or $W$ includes $B_1\setminus \{s_1\}$. Since these sets both have
cardinality $d-2$, and $|W|\le d-2$, it follows that $W$ is one of $A_1\cup C_1\setminus \{t_1\},B_1\setminus \{s_1\}$.
If $W=A_1\cup C_1\setminus \{t_1\}$, we may assume that $t_1\in P$; and since $B_1\cap W=\emptyset$, and all vertices in $B_1$
except $b_1$ are adjacent to $t_1$, it follows that $B_1\setminus \{b_1\}\subseteq P$; and in particular $u_1\in P$, and since $u_1,s_1$
are adjacent it follows that $s_1\in P$, contradicting that $Q\cap (A_1\cup B_1\cup C_1)\ne \emptyset$. Thus $W=B_1\setminus \{s_1\}$.
We may assume that $s_1\in P$, and so $A_1\cup C_1\setminus \{t_1\}\subseteq P$; and since $Q\cap (A_1\cup B_1\cup C_1)\ne \emptyset$,
it follows that $t_1\in Q$. Since $G[B_1\cup C_1'\cup \{t_1\}]$ is connected and none of its vertices are in $W$, it follows
that $B_1\cup C_1'\subseteq Q$; but there is a vertex in $C_1\setminus \{t_1\}$, and it belongs to $P$ and has a neighbour in $C_1'$,
a contradiction. This proves (1).
\\
\\
(2) {\em There do not exist $i,j\in \{1\ll m\}$ such that $V(R_i)\subseteq P\cup W$ and $V(R_j)\subseteq Q\cup W$.}
\\
\\
Suppose such $i,j$ exist; then $i\ne j$, since $|V(R_i)|>|W|$. Now the graph $D$ is $d$-connected, and so there are $d$ paths of
$G_D$, pairwise vertex-disjoint, between $V(R_i)$ and $V(R_j)$ (note that these paths are vertex-disjoint and not just 
internally disjoint, since the edges joining $V(R_i)$ to $V(G_D)\setminus V(R_i)$ form a matching, and the same for $R_j$).
But then one of these paths is disjoint from $W$, and so its vertex set is a subset of $P$ or a subset of $Q$, in either
case a contradiction since $V(R_i)\subseteq P\cup W$ and $V(R_j)\subseteq Q\cup W$. This proves (2).
\\
\\
(3) {\em For $1\le i\le m$, if $P,Q$ both have nonempty intersection with $V(R_i)$ then 
one of $A_1, A_i'$ is a subset of $P\cup W$ and the other a subset of $Q\cup W$,
and $|W\cap (C_i\cup C_i')|\ge d/2-1$.}
\\
\\
By (1) we may assume that $A_i\cup B_i\cup C_i\subseteq P\cup W$ and $A_i'\cup B_i'\cup C_i'\subseteq Q\cup W$, and the first
claim follows; and since there
is a matching between $C_i, C_i'$ of cardinality $d/2-1$, it follows that $|W\cap (C_i\cup C_i')|\ge d/2-1$. This proves (3).

\bigskip
From (2) we may assume that $P\cap V(R_i)\ne \emptyset$ for $1\le i\le m$; and so from (3), there are at most two
values of $i\in \{1\ll m\}$ such that $Q\cap V(R_i)\ne \emptyset$. Since $Q\ne \emptyset$, we may assume that
$Q\cap V(R_1)\ne \emptyset$, and from (3), $A_1\subseteq Q\cup W$ and $|W\cap (C_1\cup C_1')|\ge d/2-1$. 
If $Q\cap V(R_2)\ne \emptyset$, then similarly $|W\cap (C_2\cup C_2')|\ge d/2-1$, and so 
$W\subseteq C_1\cup C_1'\cup C_2\cup C_2'$; but some vertex in $A_1$ has a neighbour in $V(R_i)$ where $i\ne 1,2$,
and this provides an edge between $Q,P$, a contradiction. So 
$V(R_i)\subseteq P\cup W$ for $2\le i\le m$. There is a matching of cardinality $d/2$ between $A_1$ and $V(G_D)\setminus V(R_1)$,
and one of these edges has no end in $W$, and therefore joins $Q,P$, a contradiction. This completes the proof that $G_D$
is $(d-1)$-connected.

We observe that if $D,D'$ are nonisomorphic graphs, both $d$-regular and $d$-connected, then $G_D$ and $G_{D'}$ are
nonisomorphic; because
with notation as before, every subgraph of $G_D$ isomorphic to $R$ is one of the graphs $R_1\ll R_m$, and so $D$
can be obtained from $G_D$ by contracting all edges of every $R$-subgraph of $G_D$.

Now let $c>0$, and let $b=c^{4d-4}$. By \ref{manygraphs}, for infinitely many $m>0$
there are more than $b^m$ graphs $D$ on $m$ vertices that are $d$-regular and $d$-connected; and
so all the corresponding graphs $G_D$ are distinct. Hence, since each such $G_D$ has $(4d-4)m=n$ vertices,
for infinitely many $n>0$ there are more than $b^m=c^n$ $d$-regular $(d-1)$-connected graphs with $n$ vertices,
not containing a subdivision of $K_{d+1}$. This completes the proof when $d$ is even, and so proves~\ref{gadget}.~\bbox

This result \ref{gadget} is about subdivisions of $K_{d+1}$, but we believe we have a proof that
the analogous statement is true for subdivisions of $H$,
where $H$ is any $d$-regular graph. We omit the proof, which is similar but more complicated.

\section{Tangles}

A {\em separation} of a graph $G$ is a pair $(A,B)$ of subgraphs with union $G$ and with $E(A\cap B)=\emptyset$; and its {\em order}
is $|V(A\cap B)|$. Let $\theta\ge 1$ be an integer. A {\em tangle} in a graph $G$ of {\em order $\theta$} is a set $\mathcal{T}$ of separations of $G$, all of order
less than $\theta$, such that:
\begin{itemize}
\item
for every separation $(A, B)$ of order $< \theta$, one of $(A,B), (B,A)$ belongs to $\mathcal T$
\item
if $(A_1, B_1), (A_2, B_2), (A_3, B_3) \in \mathcal T$ then $A_1 \cup A_2 \cup A_3\not= G$, and
\item
if $(A, B)\in \mathcal T$ then $V(A)\not= V(G)$.
\end{itemize}
These are called the ``tangle axioms''.
Tangles were central to the ``Graph Minors'' series of papers by Robertson and the third author, and much of the theory behind them
was developed in~\cite{GM10}.

If $(A_1,B_1)\ll (A_k,B_k)$ are separations of $G$, then so is
$$(A_1\cup\cdots\cup A_k, B_1\cap \cdots\cap B_k),$$
and we call this the {\em union} of $(A_1,B_1)\ll (A_k,B_k)$. 

The following is
a consequence of theorem 2.9 of~\cite{GM10}.
\begin{thm}\label{expand}
Let $\mathcal{T}$ be a tangle of order $\theta$ in a graph $G$, and let
$(A_1, B_1)\in \mathcal{T}$.  Let $(A_2,B_2)$ be a separation of $G$ of order less than $\theta$. If either
$V(B_1)\subseteq V(B_2)$, or $V(A_2)\subseteq V(A_1)$, then $(A_2,B_2)\in \mathcal{T}$.
\end{thm}

The next result 
is theorem (8.5) of~\cite{GM10}.

\begin{thm}\label{tangledel}
Let $\mathcal{T}$ be a tangle in $G$ of order $\theta$, and let $W\subseteq V(G)$ with $|W|<\theta$. Let
$\mathcal{T}'$ be the set of all separations $(A',B')$ of $G\setminus W$ of order less than $\theta-|W|$,
such that there exists $(A,B)\in \mathcal{T}$ with $W\subseteq V(A\cap B)$ and $A'=A\setminus W$ and $B'=B\setminus W$.
Then $\mathcal{T}'$ is a tangle in $G\setminus W$ of order $\theta-|W|$.
\end{thm}
We denote the tangle $\mathcal{T'}$ in \ref{tangledel} by $\mathcal{T}\setminus W$.

\begin{thm}\label{cross}
Let $\mathcal{T}$ be a tangle in $G$ of order $\theta$, and let $(A_1,B_1), (A_2,B_2)\in \mathcal{T}$, with order $d_1,d_2$
respectively.
Let $(A,B)=(A_1\cup A_2,B_1\cap B_2)$ and $(A',B') = (A_1\cap A_2, B_1\cup B_2)$, with orders $d,d'$ respectively. Then
\begin{itemize}
\item $d+d'=d_1+d_2$;
\item if $d<\theta$ then $(A,B)\in \mathcal{T}$;
\item if $d'<\theta$ then $(A',B')\in \mathcal{T}$.
\end{itemize}
\end{thm}
\Proof
The first statement is clear.
The second statement follows from \ref{expand}, and the third from theorem 2.2 of \cite{GM10}.
This proves \ref{cross}.~\bbox

If $\mathcal{T}$ is a tangle of order $\theta$ in $G$, we say 
the {\em rank} of $X\subseteq V(G)$ relative to $\mathcal{T}$ 
is the minimum order of a separation $(A,B)\in \mathcal{T}$ with $X\subseteq V(A)$, if there is such a separation,
and $\theta$ otherwise. A set is {\em free} (relative to the tangle) if its rank equals its cardinality.
The rank of a subgraph is the rank of its vertex set. If $v\in V(G)$, a {\em $v$-nexus} in $G$ is a set 
$\mathcal{P}$ of paths of $G$, all with one end $v$. If $\mathcal{P}$ is a $v$-nexus, we write $V(\mathcal{P})$ for 
$\cup_{P\in \mathcal{P}}P$.

The next result is the main theorem of this section.

\begin{thm}\label{starrank}
Let $d',k\ge 0$ be integers. 
Let $G$ be a graph with maximum degree at most $d'$, with a tangle $\mathcal{T}$ of order at least $(k+1)^2$. Let 
$v\in V(G)$, and let $\mathcal{P}$
be a $v$-nexus, such that each member of $\mathcal{P}$ has rank at most $k$. Then $V(\mathcal{P})$ has rank at most $(kd')^k$.
\end{thm}
\Proof If $d'=0$ then $G$ has no edges, so $k=0$ and $\mathcal{P}$ has at most one path, and its rank is zero. If $d'=1$
then at most one member of $\mathcal{P}$ has maximal vertex set, and since $k\le (kd')^k$ the claim holds. Thus we may assume
that $d'\ge 2$.

Define $r_0=0$, and $r_k=(kd')^k-1$ if $k>0$. For inductive purposes we prove a slightly stronger statement, that
$V(\mathcal{P})$ has rank at most $r_k$.
We proceed by induction on $k$, and so we may assume the result holds for all $k'<k$ with $k'\ge 0$.
Let $G,\mathcal{T}, v$ and $\mathcal{P}$
be as in the theorem. For each $P\in \mathcal{P}$, since the rank of $P$ is at most $k<(k+1)^2$, there exists $(A_P,B_P)\in \mathcal{T}$
with $V(P)\subseteq V(A_P)$. 
Choose a separation $(A,B)\in \mathcal{T}$ of order at most $k$, such that the number of $P\in \mathcal{P}$ with $V(P)\subseteq V(A)$
is maximum.  We may assume that this number is at least one, and so $v\in V(A)$.

Suppose that $k=0$, and that there exists $P\in \mathcal{P}$ with $V(P)\not\subseteq V(A)$. 
Since $(A_P,B_P), (A,B)\in \mathcal{T}$, \ref{cross} implies that
$(A_P\cup A, B_P\cap B)\in \mathcal{T}$ and has order zero, contrary to the choice of $(A,B)$. Thus if $k=0$
than $V(P)\subseteq V(A)$ for each $P\in \mathcal{P}$, and so $V(\mathcal{P})$ has rank zero. Hence we may assume that $k>0$.
\\
\\
(1) {\em If $P\in \mathcal{P}$ and $V(P)\not\subseteq V(A)$ then $V(A_P\cap B_P\cap A)\ne \emptyset$.}
\\
\\
Suppose that $V(A_P\cap B_P\cap A)=\emptyset$.
From the choice of $(A,B)$, since $V(P)\subseteq V(A_P)$
and $V(P)\not\subseteq V(A)$, it follows that there exists $P'\in \mathcal{P}$ with $V(P')\subseteq V(A)$ and
$V(P')\not\subseteq V(A_P)$. Since $V(P')\subseteq V(A)$, and $V(A_P\cap B_P\cap A)=\emptyset$, it follows that
$V(A_P\cap B_P\cap P')=\emptyset$. But $v\in V(P)\subseteq V(A_P)$, and since $v\in V(P')$ and $P'$ is connected, it follows
that $P'\subseteq A_P$, a contradiction. This proves (1).

\bigskip

Let $X=V(A\cap B)$. Let $W$ be the set of vertices in $V(G)\setminus X$ with a neighbour in $X$. Thus $|W|\le kd'$.
Let $\mathcal{Q}$ be the set of all paths $Q$ of $B\setminus X$ such that $Q$ is a component of $P\setminus X$ for some
$P\in \mathcal{P}$.
Thus each member $Q$ of $\mathcal{Q}$ is a path of $G\setminus X$ and has an end in $W$.
We claim that 
\\
\\
(2) {\em Each $Q\in \mathcal{Q}$ has rank at most $k-1$ relative to $\mathcal{T}\setminus X$.}
\\
\\
First, we observe that
$(A\setminus X, B\setminus X)\in \mathcal{T}\setminus X$, of order zero. Choose $P\in \mathcal{P}$ such that $Q$ is a component of
$P\setminus X$. Consequently $V(P)\not\subseteq V(A)$. Since $(A_P,B_P)$ is a separation of order at most $k$,
it follows that there is a separation $(A_P^+, B_P^+)$ with $X\subseteq V(A_P^+\cap B_P^+)$
and $A_P^+\setminus X=A_P\setminus X$
and $B_P^+\setminus X=B_P\setminus X$, of order at most $k+|X|$; and since $(A_P,B_P) \in \mathcal{T}$, \ref{expand} implies that
$(A_P^+, B_P^+)\in \mathcal{T}$, because $k+|X|\le 2k<(k+1)^2$. Consequently $(A_P\setminus X,B_P\setminus X)\in \mathcal{T}\setminus X$.
Define $A'= (A_P\setminus X)\cup (A\setminus X)$ and $B'=(B_P\setminus X)\cap (B\setminus X)$. Then
from \ref{cross},
$(A',B')\in \mathcal{T}\setminus X.$
But
$$A'\cap B'=((A_P\cup A)\cap (B_P\cap B))\setminus X\subseteq A_P\cap B_P\cap (B\setminus X),$$
and $|V(A_P\cap B_P\cap (B\setminus X))|<|V(A_P\cap B_P)|\le k$, by (1). Thus $(A',B')$ has order at most $k-1$. Since $Q\subseteq A_P$ and $V(Q)\cap X=\emptyset$
it follows that $Q\subseteq A'$; and so $Q$ has rank at most $k-1$ relative to $\mathcal{T}\setminus X$. This proves (2).

\bigskip

Now $\mathcal{Q}$ can be partitioned
into at most $kd'$ subsets, each a $w$-nexus for some $w\in W$. Since $\mathcal{T}\setminus X$ has order at least $(k+1)^2-|X|\ge k^2$,
the inductive hypothesis implies that each such $w$-nexus
has rank at most $r_{k-1}$ relative to $\mathcal{T}\setminus X$; and so $V(\mathcal{Q})$ has rank at most $kd'r_{k-1}$
relative to $\mathcal{T}\setminus X$, and hence $V(\mathcal{Q})\cup A$ has rank at most $kd'r_{k-1}+k$ relative to $\mathcal{T}$.
Consequently $V(\mathcal{P})$ has rank at most $kd'r_{k-1}+k$ relative to $\mathcal{T}$. Since 
$kd'r_{k-1}+k\le r_k$ (because $k>0$ and $d'\ge 2$), this proves \ref{starrank}.~\bbox

\section{From minors to subdivisions}

In this section we use \ref{starrank} to prove \ref{subtominor} and \ref{4subtominor}.
We need the following, a case of theorem 7.2 of ~\cite{GM23}:

\begin{thm}\label{spider}
Let $\mathcal{T}$ be a tangle in a graph $G$, and let $W \subseteq V(G)$ be free relative to $\mathcal{T}$, with
$|W|\le w$.
Let $h \ge 1$ be an integer, and let $\mathcal{T}$ have order at least $(w+h)^{h+1} +h$.
Then there exists $W'\subseteq V(G)$ with $W\subseteq W'$ and $|W'|\le (w+h)^{h+1}$ such that for every
$(C,D)\in \mathcal{T}$ of order $<|W|+h$ with $W\subseteq V(C)$, there exists $(A',B')\in \mathcal{T}$ with
$W'\subseteq V(A'\cap B')$, such that
$|V(A'\cap B')\setminus W'| < h$ and $C \subseteq A'$.
\end{thm}

We also need:
\begin{thm}\label{maxcut}
If $\mathcal{T}$ is a tangle in a graph $G$, and $W\subseteq V(G)$ is free relative to $\mathcal{T}$, there exists
$(A_1,B_1)\in \mathcal{T}$ of order $|W|$, with $W\subseteq V(A_1)$, such that $A\subseteq A_1$ and $B_1\subseteq B$
for every $(A,B)\in \mathcal{T}$ of order $|W|$ with $W\subseteq V(A)$.
\end{thm}
\Proof Let $\mathcal{S}$ be the set of all members of $\mathcal{T}$ of order $|W|$ with $W\subseteq V(A)$.
Now $\mathcal{S}\ne \emptyset$, because $(A,B)\in \mathcal{S}$ where
$V(A)=W, E(A)=\emptyset$, and $B=G$. If $(A_1,B_1), (A_2,B_2)\in \mathcal{S}$, then their intersection has order at least $|W|$
(because otherwise it would belong to $\mathcal{T}$, by \ref{cross}, contradicting that $W$ is free); and so by \ref{cross}, their union,
$(A,B)$ say, has order at most $|W|$. Hence $(A,B)\in \mathcal{T}$ by \ref{cross}, and so it has order exactly $|W|$,
since $W$ is free; and so $(A,B)\in \mathcal{S}$. This proves that the union of every two members of
$\mathcal{S}$ is also a member of $\mathcal{S}$; and so the union of all members of $\mathcal{S}$ is a member of $\mathcal{S}$.
This proves \ref{maxcut}.~\bbox

We deduce:
\begin{thm}\label{getnodes}
Let $d,d'\ge 0$ be integers, and let $G$ be a graph with maximum degree at most $d'$.
Suppose that either $G$ is $d$-connected, or $d\in \{4,5\}$ and $G$ is $(d-1)$-connected and $d$-edge-connected.
Let $s\ge 0$ be an integer, and let $\mathcal{T}$ be a tangle in $G$ of order
$$\theta\ge (s-1)(d+1)+ (dd')^d(d'(s-1+(sd)^{d+1})+d(s-1))+(sd)^{d+1}.$$
Then there exist
distinct $z_1\ll z_s\in V(G)$, and pairwise disjoint subsets $W_1\ll W_s$ of $V(G)\setminus \{z_1\ll z_s\}$, such that
\begin{itemize}
\item for $1\le i\le s$, $|W_i|= d$, and $z_i$ is adjacent to each vertex in $W_i$; and
\item $W_1\cup\cdots\cup W_s$ is free relative to $\mathcal{T}\setminus \{z_1\ll z_s\}$ in $G\setminus \{z_1\ll z_s\}$.
\end{itemize}
\end{thm}
\Proof
We proceed by induction on $s$, and so we may assume that there exists
$Z=\{z_1\ll z_{s-1}\}\subseteq  V(G)$, and pairwise disjoint subsets $W_1\ll W_{s-1}$ of $V(G)\setminus Z$, such that
\begin{itemize}
\item for $1\le i\le s-1$, $|W_i|= d$, and $z_i$ is adjacent to each vertex in $W_i$;
\item $W_1\cup\cdots\cup W_{s-1}$ is free relative to $\mathcal{T}\setminus Z$ in $G\setminus Z$.
\end{itemize}
Let $G'=G\setminus Z$ and $\mathcal{T}'=\mathcal{T}\setminus Z$; so $\mathcal{T}'$ is a tangle in $G'$ of order $\theta-s+1$.
Let $W=W_1\cup\cdots\cup W_{s-1}$; then $|W|=(s-1)d$ and $W$
is free relative to $\mathcal{T}'$. From \ref{maxcut}, there exists $(A_1,B_1)\in \mathcal{T}'$, of order $|W|$ and
with $W\subseteq V(A_1)$, such that
$A\subseteq A_1$ and $B_1\subseteq B$ for every $(A,B)\in \mathcal{T}'$ of order $|W|$ with $W\subseteq V(A)$.
\\
\\
(1) {\em $W$ is free relative to $\mathcal{T}'\setminus \{v\}$, for each $v\in V(B_1)\setminus V(A_1)$.}
\\
\\
Suppose not; then there is a separation $(A,B)$ of $\mathcal{T}'\setminus \{v\}$ of order $<|W|$, with $W\subseteq V(A)$.
Hence there is a separation $(A',B')$ of $\mathcal{T}'$ of order $\le |W|$ with $W\subseteq V(A')$ and $v\in V(A'\cap B')$.
But $(A',B')$ has order exactly $|W|$, since $W$ is free relative to $\mathcal{T}'$; and so $A'\subseteq A_1$, 
from the property of $(A_1,B_1)$.
This contradicts that $v\in V(B_1)\setminus V(A_1)$, and so proves (1).

\bigskip

Let $w'=(sd)^{d+1}$.
By \ref{spider} applied to $G',\mathcal{T}'$ (taking $w=(s-1)d$ and $h=d$)
there exists $W'\subseteq V(G')$ with $W\subseteq W'$ and $|W'|\le w'$, such that for every
$(C,D)\in \mathcal{T}'$ of order $<sd$ with $W\subseteq V(C)$, there exists $(A',B')\in \mathcal{T}'$ with
$W'\subseteq V(A'\cap B')$, such that
$|V(A'\cap B')\setminus W'| < d$ and $C\subseteq A'$.

Let $G''=G'\setminus W'$, and $\mathcal{T}''=\mathcal{T}'\setminus W'$. Hence $\mathcal{T}''$ is a tangle in $G''$ of 
order at least $\theta-s+1-w'$.
Let $N$ be the set of vertices of $G$ that are not in $W'\cup Z$ but have a neighbour in $W'\cup Z$. Hence $N\subseteq V(G'')$. 
Since every vertex of $G$
has degree at most $d'$, it follows that $|N|\le d'|W'\cup Z|\le d'(s-1+w')$. 
Let $N'=N\cup V(A_1\cap B_1\cap G'')$. Thus $|N'|\le d'(s-1+w')+d(s-1)$.

By \ref{starrank}, for each $n\in N'$ there is a separation $(A_n,B_n)\in \mathcal{T}''$ of order at most $(dd')^d$, such that 
$P\subseteq A_n$ for each path $P$ of $G''$ of rank at most $d$ (relative to $\mathcal{T}''$) with one end $n$. Let $(A_0,B_0)$ be the union of these
separations. Thus $(A_0,B_0)$ has order at most
$(dd')^d(d'(s-1+w')+d(s-1))$, and so belongs to $\mathcal{T}''$ by \ref{cross}.
\\
\\
(2) {\em There exists $v\in V(G'')$ such that $v\notin V(A_0\cup A_1)$.}
\\
\\
Since $(A_0,B_0)\in \mathcal{T}''$, there is a separation $(A_2,B_2)$ of $\mathcal{T}'$ with $W'\subseteq V(A_2\cap B_2)$,
such that $A_2\setminus W'=A_0$ and
$B_2\setminus W'=B_0$. Its order is at most $(dd')^d(d'(s-1+w')+d(s-1))+w'$; and so the union $(A,B)$ of $(A_1,B_1), (A_2,B_2)$
has order at most
$$(dd')^d(d'(s-1+w')+d(s-1))+w'+ (s-1)d$$
and so belongs to $\mathcal{T}'$.
From the third tangle axiom, applied to $\mathcal{T}'$ and $(A,B)$,
there is a vertex $v\in V(B)\setminus V(A)$. In particular, $v\notin W'$ since $W'\subseteq V(A_2)$, and so $v\in V(G'')$.
This proves (2).
\\
\\
(3) {\em There exists $u\in V(G'')\setminus V(A_1)$ such that there
is no $(A,B)\in \mathcal{T}''$ of order $<d$ with $u\in V(A)\setminus V(B)$.}
\\
\\
Choose $v$ as in (2); we may therefore assume that there is a separation $(A,B)\in \mathcal{T}''$ of order $<d$ with $v\in V(A)\setminus V(B)$.
Choose $(A,B)$ with $B$ minimal.
Let $C$ be the component of $A$ that contains $v$, and suppose first that $V(C)\cap N'\ne \emptyset$.
Hence there is a path of $A$ between $v$ and $N'$, say $P$.
Thus $P$ has rank at most $d-1$ relative to $\mathcal{T}''$, 
since $P\subseteq A$; and so $P\subseteq A_0$, contradicting that $v\notin V(A_0)$. This proves that $V(C)\cap N'=\emptyset$,
and in particular, $G$ is not $d$-connected. 

Consequently $d\in \{4,5\}$, and $G$ is $(d-1)$-connected and $d$-edge-connected. Thus $(A,B)$ has order $d-1$. 
Since $G$ is 
$d$-edge-connected, there are at least $d$ edges of $G$ between $V(C)\setminus V(B)$ and its complement in $V(G)$. 
Since all of these edges belong to $G''$
(since $V(C)\cap N=\emptyset$), they are all between $V(C)\setminus V(B)$ and $V(A\cap B)$; and so some vertex $u\in V(A\cap B\cap C)$ 
has at least two neighbours $u_1,u_2\in V(C)\setminus V(B)$. In particular, $u\in V(C)$; let $P$ be a path of $C$
between $u,v$. If $u\in V(A_1)$, then since $v\notin V(A_1)$,
some vertex of $P$ belongs to $V(A_1\cap B_1)$, and since $P$ has rank at most $d-1$ relative to $\mathcal{T}''$, 
it follows that $P\subseteq A_0$ from the definition of $(A_0,B_0)$, contradicting that $v\notin V(A_0)$. Thus $u\notin V(A_1)$.

Suppose that there is a separation $(A',B')$ of $G''$ of order $<d$ with $u\in V(A')\setminus V(B')$. We may assume that
$A'$ is connected (for instance, by choosing
$(A',B')$ with $A'$ minimal). If $N'\cap V(A')\ne \emptyset$, then there is a path of $A'$ between $N'$ and $u$, and since there is 
a path of $C$ between $u,v$ included in $A$, it follows that there is a path $P$ between $N', v$ included in $A\cup A'$. 
But $P$ has rank at most the order of $(A\cup A', B\cap B')$, and hence at most $2d-5\le d$, and so $P\subseteq A_0$
from the definition of $(A_0,B_0)$, contradicting that $v\notin V(A_0)$. This proves that  $N'\cap V(A')= \emptyset$.
If $(A\cup A',B\cap B')$  has order at most $d-1$, then it belongs to $\mathcal{T}''$ by \ref{cross}, contradicting the minimality
of $B$ since $A\cup A'\ne A$. Thus $(A\cup A',B\cap B')$ has order at least $d$, and since $(A,B),(A',B')$ have order at most $d-1$,
it follows that $(A\cap A',B\cup B')$ has order at most $d-2$. Since $G$ is $(d-1)$-connected and $N'\cap V(A')= \emptyset$, it follows that 
$V(A\cap A')\subseteq V(B\cup B')$, and in particular $u_1,u_2\in V(B')$. Thus $|V(A'\cap B')\setminus V(B)|\ge 2$. Since 
$|V(A'\cap B')|\le d-1$, it follows that $|V(A'\cap B')\cap V(B)|\le d-3$. Since $(A\cap A', B\cap B')$ has order
greater than the order of $(A,B)$, it follows that 
$$2\ge d-3\ge |V(A'\cap B')\setminus V(A)|>|V(A\cap B)\setminus V(B')|\ge 1,$$
and so equality holds throughout; and in particular, $d=5$, $V(A\cap B)\setminus V(B')=\{u\}$, and $|V(A'\cap B')\setminus V(A)|=2$,
and so $V(A\cap B\cap A'\cap B')=\emptyset$. Thus the separation $(A'\cap B, A\cup B')$ has order at most three, and since 
$G$ is $4$-connected and $N'\cap V(A')= \emptyset$, it follows that $V(A'\cap B)\subseteq V(B'\cup A)$, and in particular, all neighbours of $u$ belong
to $V(A'\cap B')$; which is impossible since $G$ is $d$-edge-connected. Thus there is no such separation $(A',B')$.
This proves (3).

\bigskip
Define $z_s=u$, where $u$ is chosen as in (3).
\\
\\
(4) {\em There is no $(C,D)\in \mathcal{T}'$ of order $<sd$ such that $W\subseteq V(C)$ and $z_s\in V(C)\setminus V(D)$.}
\\
\\
For suppose that there is such a separation $(C,D)$. From the choice of $W'$,
there exists $(A',B')\in \mathcal{T}'$ with $W'\subseteq V(A'\cap B')$
such that
$|V(A'\cap B')\setminus W'| < d$ and $C\subseteq A'$. Hence $z_s\in V(A')$. Choose $(A',B')$ with $B'$ minimal, and suppose that
$z_s\in V(B')$. From the minimality of $B'$, it follows that $(A',B'\setminus \{z_s\})\notin \mathcal{T}'$, and so by \ref{expand},
$z_s$ is adjacent in $G'$ to some vertex $b\in V(B')\setminus V(A')$. But $z_s\in V(C)\setminus V(D)$, and so $b\in V(C)\subseteq V(A')$,
a contradiction. Thus $z_s\notin V(B')$, and so $z_s\in V(A')\setminus V(B')$. But $(A'\setminus W', B'\setminus W')\in \mathcal{T}''$,
contrary to (3). This proves (4).

\bigskip
Let $P$ be the set of all neighbours of $z_s$ in $V(G)\setminus (Z\cup W)$.
\\
\\
(5) {\em There is no $(A,B)\in \mathcal{T}'\setminus \{z_s\}$ of order $<sd$ such that $W\cup P\subseteq V(A)$.}
\\
\\
For suppose there exists such a separation $(A,B)$ of $G'\setminus \{z_s\}$. Hence there is a separation $(C,D)\in \mathcal{T}'$
with $z_s\in V(C\cap D)$, such that $C\setminus \{z_s\} = A$ and $D\setminus \{z_s\}=B$. Choose $(C,D)$ with $C$ maximal. By \ref{expand},
every edge of $G'$ incident with $z_s$ and with its other end in $V(C)$ belongs to $C$.
But every neighbour of $z_s$ in $G'$ belongs to $P\cup W$ and hence to $V(C)$;
and so no edge of $D$ is incident with $z_s$. Consequently $(C, D\setminus \{z_s\})$ is also a separation of $G'$, and by \ref{expand},
it belongs to $\mathcal{T}'$. Its order is that of $(A,B)$ and hence less than $sd$, contrary to (4). This proves (5).

\bigskip

By theorem 12.2 of~\cite{GM10}, the subsets of $V(G'\setminus \{z_s\})$ that are free relative to $\mathcal{T}'\setminus \{z_s\}$ are the independent sets of a
matroid, $M$ say, of rank $\theta-s$ (since $\mathcal{T}'\setminus \{z_s\}$ has order $\theta-s$); and by the same theorem,
the set $P\cup W$ has rank
at least $sd$ in $M$ because of (5). Since $W$ is free relative to $\mathcal{T}'\setminus \{z_s\}$ by (1), and hence independent in $M$,
it can be extended to an independent subset of $W\cup P$ of cardinality $sd$; and so
there exists $W_s\subseteq P$ of cardinality $d$ such that $W\cup W_s$ is free relative to $\mathcal{T}'\setminus \{z_s\}$.
This completes the inductive argument, and so proves \ref{getnodes}.~\bbox

We will need theorem 5.4 of \cite{GM13}, which states:
\begin{thm}\label{useclique}
Let $G$ be a graph and let $Z\subseteq V(G)$. Let $k\ge \lfloor \frac32 |Z|\rfloor$ and let 
$C_1\ll C_k$ be connected subgraphs of $G$, mutually vertex-disjoint, such that for $1\le i<j\le k$ there is
an edge between $C_i$ and $C_j$. Suppose that there is no separation $(A,B)$ of $G$ of order $<|Z|$ with $Z\subseteq V(A)$
and $A\cap C_i$ null for some $i\in \{1\ll k\}$. Then for every partition $(Z_i:1\le i\le n)$ of $Z$ into nonempty subsets,
there are $n$ connected subgraphs $T_1\ll T_n$ of $G$, mutually vertex-disjoint and with $V(T_i)\cap Z=Z_i$ for $1\le i\le n$.
\end{thm}

Finally, we deduce \ref{subtominor} and \ref{4subtominor}, which we combine in:
\begin{thm}\label{subtominor2}
Let $d,d'\ge 0$, and let $H$ be a graph with maximum degree at most $d$.
Then there exists $t\ge 0$ such that for every graph $G$, if 
\begin{itemize}
\item $G$ has maximum degree at most $d'$, 
\item $G$ contains $K_t$ as a minor, and 
\item either $G$ is $d$-connected,
or $d\in \{4,5\}$ and $G$ is $(d-1)$-connected and $d$-edge-connected, 
\end{itemize}
then $G$ contains a subdivision of $H$.
\end{thm}
\Proof
Let $s=|V(H)|$, and we assume $V(H)=\{h_1\ll h_s\}$. Let 
$$\theta\ge (s-1)(d+1)+ (dd')^d(d'(s-1+(sd)^{d+1})+d(s-1))+(sd)^{d+1},$$
let $k=\lceil 3\theta/2\rceil$, and let $t=k+s$. Now let $G$ be as in the theorem.
Choose $t$ disjoint connected subgraphs $C_1\ll C_t$ such that for $1\le i<j\le t$ there is an edge
of $G$ between $V(C_i)$ and $V(C_j)$. 

For each separation $(A,B)$ of $G$ of order $<\theta$, exactly one of $A\setminus V(B),B\setminus V(A)$ includes one of $C_1\ll C_t$;
let $\mathcal{T}$ be the set of all such $(A,B)$ where $B\setminus V(A)$ includes one of $C_1\ll C_t$. Then $\mathcal{T}$
is a tangle in $G$ of order $\theta$, for instance by the argument of theorem 4.4 of~\cite{GM10}.
By \ref{getnodes} applied to $\mathcal{T}$, there exist $z_1\ll z_s\in V(G)$, and pairwise disjoint subsets $W_1\ll W_s$ of $V(G)\setminus \{z_1\ll z_s\}$, 
such that
\begin{itemize}
\item for $1\le i\le s$, $|W_i|= d$, and $z_i$ is adjacent to each vertex in $W_i$; and
\item $W_1\cup\cdots\cup W_s$ is free relative to $\mathcal{T}\setminus \{z_1\ll z_s\}$ in $G\setminus \{z_1\ll z_s\}$.
\end{itemize}
For $1\le i\le s$ and each edge $e$ of $H$ incident with $h_i$, choose $y_{i,e}\in W_i$, in such a way that all the vertices 
$y_{i,e}$ are distinct (this is possible since $|W_i|=d$ and $h_i$ has degree at most $d$ in $H$). 

Let $Z=W_1\cup\cdots\cup W_s$, and for each edge $e=h_ih_j$ of $H$, let $Z_e=\{y_{i,e},y_{j,e}\}$. Then $(Z_e:e\in E(H))$
is a partition of $Z$. We may assume that none of $C_1\ll C_{k}$ contain any of $z_1\ll z_s$, since $C_1\ll C_t$
are pairwise disjoint, and so $C_1\ll C_s$ are all subgraphs of $G'$, where $G'=G\setminus \{z_1\ll z_s\}$. 
Suppose that there is a separation $(A',B')$ of $G'$ of order $<|Z|$ with $Z\subseteq V(A')$
and $A'\cap C_i$ null for some $i\in \{1\ll k\}$. It follows that there is a separation $(A,B)$ of $G$ of order less than $|Z|+s$,
with $\{v_1\ll v_s\}\subseteq V(A\cap B)$, such that $A\setminus \{z_1\ll z_s\}=A'$ and $B\setminus \{z_1\ll z_s\}=B'$. Consequently
$A\cap C_i$ is null, since $z_1\ll z_s\notin V(C_i)$; and so $(A,B)\in \mathcal{T}$, since $|Z|+s=s(d+1)< \theta$. 
Hence $(A',B')\in \mathcal{T}\setminus \{z_1\ll z_s\}$, contradicting that $Z$ is free relative to 
$\mathcal{T}\setminus \{z_1\ll z_s\}$ in $G\setminus \{z_1\ll z_s\}$. Thus there is no such $(A',B')$.

From \ref{useclique} applied to $G', \mathcal{T}\setminus \{z_1\ll z_s\}$, $Z$ and the partition $(Z_e:e\in E(H))$, 
it follows that for each edge $e$ of $H$ there is a connected subgraph $P_e$ of $G'$ containing the two vertices of $Z_e$,
such that the subgraphs $P_e\;(e\in E(H))$ are pairwise vertex-disjoint. By choosing each $P_e$ minimal we may assume each $P_e$
is a path joining the two members of $Z_e$. But then adding the vertices $z_1\ll z_s$ and the edges between each $z_i$ 
and the corresponding $W_i$, to the union of these paths, gives a subdivision
of $H$. This proves \ref{subtominor2}.~\bbox

\end{document}